\documentclass[reqno,10pt,oneside]{amsart}
\usepackage{amssymb,latexsym}
\usepackage{epsf}
\def\xx{\bar{x}}
\def\yy{\bar{y}}
\def\CC{\mathbb{C}}
\def\ZZ{\mathbb{Z}}

\def\ct{\mathop{\rm CT}}
\def\CT{\mathop{\rm CT}}

\def\bous{Bousquet-M\'{e}lou}

\newtheorem{thm}{Theorem}[section]

\newtheorem{lem}[thm]{Lemma}

\allowdisplaybreaks[1]

\makeatletter\@addtoreset{equation}{section}\makeatother

\begin{document}

\title{Proof of a Conjecture on the Slit Plane Problem}
\author{Guoce Xin}
\address{Department
of Mathematics\\
Brandeis University\\
Waltham MA
02454-9110}
\email{maxima@brandeis.edu}
\date{April 14, 2003}
\begin{abstract}
Let $a_{i,j}(n)$ denote the number of walks in
$n$ steps from $(0,0)$ to $(i,j)$, with steps $(\pm 1,0)$ and
$(0,\pm 1)$, never touching  a point $(-k,0)$ with $k\ge 0$ after
 the starting point.
\bous\ and Schaeffer conjectured   a closed form for
the number $a_{-i,i}(2n)$ when $i\ge 1$.
 In this paper, we prove their conjecture, and
give a formula for $a_{-i,i}(2n)$ for $i\le -1$.
\end{abstract}

\maketitle

\section{Introduction and Theorems}
The problem of { walks on the slit plane} was first studied by
M. \bous\ and G. Schaeffer in \cite{bous}. See also \cite{bouso}.

Let $a_{i,j}(n)$ denote the number of walks in
$n$ steps from $(0,0)$ to $(i,j)$, with steps $(\pm 1,0)$ and
$(0,\pm 1)$, never touching  a point $(-k,0)$ with $k\ge 0$ after
 the starting point. These are called {\em walks on the slit plane}.

Let $\xx$ denote $x^{-1}$ and $\yy$ denote $y^{-1}$.
In \cite{bous}, the authors showed (Theorem $1$) that
\begin{align}
\label{e-th1}
S(x,y;t)&=\sum_{n\ge 0}\sum_{i,j\in \ZZ}
a_{i,j}(n)x^iy^jt^n=\frac{(1-2t(1+\xx)+\sqrt{1-4t})^{1/ 2} (1+2t(1-\xx)+\sqrt{1+4t})^{1/ 2}}
{2(1-t(x+\xx+y+\yy))},
\end{align}
where $S(x,y;t)$ is the complete generating function for walks on
the slit plane.

The authors also conjectured a closed form for $a_{-i,i}(2n)$ for $i\ge 1$.
By reflecting in the $x$-axis, we see that $a_{-i,i}(2n)=a_{-i,-i}(2n)$, the closed
form of which is given as \eqref{conj1}
in the following theorem.
\begin{thm}\label{t-conj}
For $i\ge 1 $ and $n\ge i$, we have
\begin{align}
a_{-i,-i}(2n) & =\frac{i}{2n}{2i\choose i}{n+i\choose 2i}\frac{{4n\choose 2n}}{
{2n+2i\choose 2i}}, \label{conj1} \\
a_{i,i}(2n) &= a_{-i,-i}+4^n\frac{i}{n}\binom{2i}{i}\binom{2n}{n-i}
  \label{conj2}.
\end{align}
\end{thm}

We will prove this
theorem in the next section.
Theorem \ref{t-lagrange1} below is a basic tool
to prove the conjecture.

There are two key steps in proving the conjecture that might be worth
mentioning: one is using Theorem \ref{t-lagrange1}
to obtain the generating function \eqref{conj3} that involves $a_{i,i}(2n)$ for all integers $i$;
the other is guessing the formula \eqref{conj2}.

Let $R$ be a commutative ring with unit, and $R[x,\xx][[t]]$
the ring of formal power series in $t$ with coefficients Laurent
polynomials in $x$. An element of $R[x,\xx][[t]]$ is written
as $f(x;t)$, to emphasize that $f(x;t)$ is regarded as
a power series in $t$.

If $f(x;t)\in R[x,\xx][[t]]$, then it can be written as
$$f(x;t)=\sum_{n\ge 0}\sum_{i\in \ZZ} f_{i,n} x^it^n.$$

Let
$\ct_x f(x;t)$ denote the constant term of $f(x;t)$ in $x$, i.e.,
$$\ct_x f(x;t)=\sum_{n\ge 0} f_{0,n} t^n.$$

Note that if we are working in $R[x,\xx][[t]]$, and if $u\in tR[[t]]$, then
$\frac{x}{x-u}$ has to be interpreted
as
$$\frac{x}{x-u}=\frac{1}{1-u/x}=\sum_{n\ge 0} u^n/x^n. $$

So
$$\ct_x \frac{x}{x-u} x^k = \left\{ \begin{array}{ll}
u^k & \text{if } k\ge 0, \\
0 & \text{if } k<0.
\end{array}\right. $$

By linearity, we have the following:
\begin{lem}\label{l-l11}
Let $Q(x;t)\in R[x][[t]]$ be a formal power series in $t$,
with coefficients  in the polynomial ring $R[x]$.
If $u=u(t)\in tR[[t]]$ is a formal power series in $t$ with constant term $0$,
then
\begin{align}\label{e-ctqq}
\ct_x \frac{x}{x-u} Q(x;t) = Q(u;t).
\end{align}
\end{lem}

The following lemma is a well-known result. See, e.g., \cite[Theorem 4.2]{gessel}.
\begin{lem}
If $G(x,t)\in R[[x,t]]$, then there is a unique $X=X(t)$ in $tR[[t]]$ such that
$X-tG(X,t)=0$.
\end{lem}

\begin{thm}\label{t-lagrange1}
Let $G(x,t), F(x,t) \in R[[x,t]]$, and let $X=X(t)$ be the unique element
in $tR[[t]]$ such that $X-tG(X,t)=0$. Then
\begin{equation}
\CT_{x} \frac{x}{x-tG(x,t)} F(x,t)= \left. \frac{F(X,t)}{1-t \displaystyle{\partial\over
\partial x} G(x,t)}\right|_{x=X}.
\end{equation}
\end{thm}
\begin{proof}
Write $G(x,t)= \sum_{n\ge 0} a_n(t) x^n$. Then
$$\frac{x-t G(x,t)}{x-X}=\frac{x-tG(x,t)-(X-tG(X,t))}{x-X}
=1-t\sum_{n\ge 0} a_n(t) (x^{n-1}+x^{n-2}X+\cdots
+X^{n-1}),
$$
which is an element in $R[[x,t]]$ with constant term $1$.
Setting $x=X$, we get
$$\left. \frac{x-t G(x,t)}{x-X}\right|_{x=X} = 1-\left. t{\partial\over \partial x} G(x,t)\right|_{x=X}
.$$

By Lemma \ref{l-l11}, we have
\begin{align*}
\CT_{x} \frac{x}{x-tG(x,t)} F(x,t) &= \CT_x \frac{x}{x-X} \left(\frac{x-tG(x,t)}{x-X}\right)^{-1}
F(x,t) \\
&=\left. \left(\frac{x-tG(x,t)}{x-X}\right)^{-1}
F(x,t) \right|_{x=X} \\
&= \left. \frac{F(X,t)}{1-t \displaystyle{\partial\over
\partial x} G(x,t)} \right|_{x=X} .
\end{align*}
\end{proof}

Theorem \ref{t-lagrange1} is a generalization of Lagrange's inversion formula.
If we set $F$ and $G$ to be independent of $t$, we can easily derive Lagrange's
inversion formula. See \cite[Theorem 5.4.2]{stanley}.
This topic will be explored further in \cite{guoce}.

\vspace{3mm}
\section{The Proof of the Conjecture}
Let $$C(t)=\sum_{n\ge 0} C_n t^n= \frac{1-\sqrt{1-4t}}{2t}$$ be the Catalan
generating function, and
let $$u=tC(t)C(-t)=\frac{\sqrt{1+4t}-1}{\sqrt{1-4t}+1}.$$

Much of the computation here involves rational functions of
$u$. We shall use the following facts from \cite{bous}.
$$\CC(u)=\CC(t,\sqrt{1-4t},\sqrt{1+4t}),$$
$$C(t)=\frac{1+u^2}{1-u}, \quad C(-t)=\frac{1+u^2}{1+u},
\quad C(4t^2)=\frac{(1+u^2)^2}{(1-u^2)^2}.$$

We shall
prove Theorem \ref{t-conj}
by computing the diagonal generating function $F(y;t)$. More precisely, let
$$F(y;t)=\sum_{n\ge 0}\sum_{i\in \ZZ} a_{i,i}(2n) y^i t^{2n}.$$

Since $S(x,y\xx;t)$ belongs to $\CC[x,y,\xx,\yy][[t]]$, it is easy to check that
\begin{align}\label{e-dd1}
F(y;t) =\CT_x S(x,y\xx;t)= \CT_x S(\xx, xy;t).
\end{align}

\begin{lem}\label{l-conj1}
\begin{equation}\label{conj3}
F(y;t)=\frac{\left[\frac{1+\sqrt{1-4t}}{2}-t\left(1+\frac{1-\sqrt{1-4t^2(1+y)^2/y}}
{2t(1+y)}\right)\right]^{1/2}
\left[\frac{1+\sqrt{1+4t}}{2}+t\left(1-\frac{1-\sqrt{1-4t^2(1+y)^2/y}}{2t(1+y)}\right)
\right]^{1/2} }{\sqrt{1-4t^2(1+y)^2/y}}.
\end{equation}
\end{lem}
\begin{proof}
Using \eqref{e-dd1} and \eqref{e-th1}, we get
\begin{align*}
F(y;t)=& \CT_x S(\xx, xy;t) \\
=&\CT_x \frac{1}
{2(1-t(x+\xx+x y+\xx \yy))}(1-2t(1+x)+
\sqrt{1-4t})^{1/2}(1+2t(1-x)+\sqrt{1+4t})^{1/2}\\
=&\CT_x \frac{x}{2(x-t(x^2+1+x^2y+\yy))} (1-2t(1+x)+\sqrt{1-4t})^{1/2}(1+2t(1-x)+\sqrt{1+4t})^{1/2}.
\end{align*}

Applying Theorem
\ref{t-lagrange1} with $R=\CC[y,\yy]$, this becomes
$$ \frac{1}{2\left(1-t (2X+2X y)\right)}
(1-2t(1+X)+\sqrt{1-4t})^{1/2}(1+2t(1-X)+\sqrt{1+4t})^{1/2},
$$
where $X=X(t)$ is the unique solution in $tR[[t]]$ such that $X=t(X^2+1+\yy+X^2y)$.
We can solve for  $X$ by the quadratic formula:
$$X= \frac{1-\sqrt{1-4t^2(1+y)^2/y}}{2t(1+y)}.$$

Equation \eqref{conj3}  then follows.
\end{proof}

It is clear that for any $G(y;t)\in R[y,\yy][[t]]$, there is a unique
decomposition $G(y;t)=G_+(y;t)+G_0(t)+G_-(\yy,;t)$, such that
$G_+(y;t),G_-(y;t)\in yR[y][[t]]$ and $G_0(t)\in R[[t]]$.

Our task now is to find this decomposition of $F(y;t)$.
 There is no general theory to do this.
For this particular $F(y;t)$, thanks to the work
of \bous\ and Schaeffer,
we can guess the formulas for $F_+$ and $F_-$ and prove them.

The variable $s$ defined by the following is useful:
\begin{equation}
s=tC(4t^2)=\frac{u}{1-u^2} \text{ and } t=\frac{s}{1+4s^2} \label{e-of-z}.
\end{equation}
Note that $s$ is also $S_{0,1}(t)$, the generating function
of walks on the slit plane that end at $(0,1)$. See \cite[p. 11]{bous}.

\begin{lem}\label{l-conj2}
We have the decomposition
$$F(y;t)=F_+(y,t)+1+F_-(\yy,t),$$
where
\begin{align}
F_+(y,t) & =F_-(y,t)+\frac{1}{2}((1-4s^2y)^{-1/2}-1),\label{e-f+} \\
F_-(y,t) &=\frac{(1-u^2)s^2yC(s^2y)}{1+u^2C^2(s^2y)s^2y}
\frac{1}{\sqrt{1-4s^2y}}\label{conj4} .
\end{align}
\end{lem}
\begin{proof}
Let
\begin{align*}
T(y;t)&=\frac{(1-u^2)s^2yC(s^2y)}{1+u^2C^2(s^2y)s^2y}
\frac{1}{\sqrt{1-4s^2y}}+\frac{1}{2}((1-4s^2y)^{-1/2}-1)+1+
\frac{(1-u^2)s^2\yy C(s^2\yy)}{1+u^2C^2(s^2\yy)s^2\yy}
\frac{1}{\sqrt{1-4s^2\yy}}.
\end{align*}

From Lemma \ref{l-conj1}, we have
\begin{align*}
F(y;t)&=\frac{\left[\frac{1+\sqrt{1-4t}}{2}-t\left(1+\frac{1-\sqrt{1-4t^2(1+y)^2/y}}
{2t(1+y)}\right)\right]^{1/2}
\left[\frac{1+\sqrt{1+4t}}{2}+t\left(1-\frac{1-\sqrt{1-4t^2(1+y)^2/y}}{2t(1+y)}\right)
\right]^{1/2} }{\sqrt{1-4t^2(1+y)^2/y}}.
\end{align*}

Therefore, it suffices to show that $T(y;t)=F(y;t)$.
Since it is easy to see that $T(y;0)=F(y;0)=1$, the proof will be
completed by showing that $T^2(y;t)-F^2(y;t)=0$.

Using the variable $u$, we can get rid of the radicals $\sqrt{1-4t}$ and
$\sqrt{1+4t}$ by the following:
$$\sqrt{1-4t}=\frac{1-2u-u^2}{1+u^2}, \text{ and } \sqrt{1+4t}=\frac{1+2u-u^2}{1+u^2}.$$

The radicals left are
$D=\sqrt{1-4s^2y}$, $E=\sqrt{1-4s^2\yy}$, and $\sqrt{1-4t^2(1+y)^2/y}$, which is easily
checked to be
equal to
$DE$.

Rewriting $T^2-F^2$ in terms of $u,D,E$, we get a rational function
of $u,D,E$. For $i=1,2$ (the degrees in $D$ and
$E$ are both $4$), replacing $D^{2i}$ by $(1-4s^2y)^i$, $D^{2i+1}$
by $(1-4s^2y)^iD$, $E^{2i}$ by $(1-4s^2\yy)^i$, and
$E^{2i+1}$ by $(1-4s^2\yy)^iE$, we
find that the expression reduces to $0$.
\end{proof}

Now we need to show the following.
\begin{lem}\label{l-conj3}
\begin{align}\label{conj7}
F_-(y,t)&=\sum_{n\ge 0} \sum_{i\ge 1} b_{i}(2n) t^ny^i ,
\end{align}
where
\begin{align}
b_i(2n)&=\frac{i}{2n}{2i\choose i}{n+i\choose 2i}\frac{{4n\choose 2n}}{
{2n+2i\choose 2i}}.
\end{align}
\end{lem}
We will give two proofs of this lemma. The first one starts from a formula in
\cite{bous}. We include it here as an example of computing the generating
function by Theorem \ref{t-lagrange1}. The second proof is self-contained, and
is simpler.

Let
\begin{align}
\label{e-dd2}
f(y,t)=\sum_{n\ge 1} \sum_{i\ge 1} b_i(2n) t^ny^i .
\end{align}

We need to show that $F_-(y,t)=f(y,t)$.
\begin{proof}[First Proof of Lemma \ref{l-conj3}.]
It was stated in \cite{bous} that
\begin{equation}\label{conj8}
\sum_{n\ge 0} b_i(2n)
t^n=\frac{(-1)^i}{(1-u^2)^{2i-1}}\sum_{k=i}^{2i-1} {2i-1\choose k} (-1)^k
u^{2k}.
\end{equation}

Let $s$ be as in \eqref{e-of-z}. Using the following fact
$${n\choose k}= \ct_\alpha \frac{1}{\alpha^k} (1+\alpha)^n,$$
we can compute $f(y,t)$ by Theorem \ref{t-lagrange1}:
\begin{align*}
f(y,t)
&=\sum_{i\ge 1} \frac{(-1)^i}{(1-u^2)^{2i-1}}\sum_{k=i}^{2i-1}
{2i-1\choose k} (-1)^k u^{2k} y^i \\
&=\sum_{i\ge 1} \frac{(1-u^2)(-1)^i}{(1-u^2)^{2i}} \sum_{r=0}^{i-1}
{2i-1\choose i+r} (-1)^{i+r} u^{2i+2r} y^i , \mbox{ where } r=k-i\\
&= (1-u^2) \sum_{r\ge 0} (-1)^ru^{2r}\sum_{i\ge r+1}
{2i-1\choose i-1-r} \frac{u^{2i}} {(1-u^2)^{2i}}y^i\\
&= (1-u^2)\sum_{r\ge 0}(-u^2)^r\sum_{i\ge r+1} \CT_\alpha
(1+\alpha )^{2i-1}\left(\frac{1}{\alpha}\right)^{i-1-r} (s^2y)^i \\
&= \CT_\alpha (1-u^2)\sum_{r\ge 0} \frac{\alpha }{1+\alpha }
(-u^2)^r\alpha ^r \sum_{i\ge r+1} \frac{(1+\alpha)^{2i}}{\alpha^i} (s^2y)^i \\
&= \CT_\alpha \frac{\alpha}{1+\alpha}(1-u^2)\sum_{r\ge 0}(-u^2\alpha)^r
\left(\frac{(1+\alpha)^2}{\alpha} s^2y\right)^{r+1}
\frac{1}{1-\displaystyle\frac{(1+\alpha)^2}{\alpha}s^2y} \\
&= \CT_\alpha (1-u^2)(1+\alpha) s^2y \frac{1}{1+u^2(1+\alpha)^2s^2y} \cdot
\frac{1}{1-\displaystyle\frac{(1+\alpha)^2}{\alpha}s^2y} .
\end{align*}

Now
$$(1-u^2)(1+\alpha) s^2y \displaystyle\frac{1}{1+u^2(1+\alpha)^2s^2y}
$$ is a power series in $t$ with coefficients in $\CC[y][\alpha]$,
and
\begin{align*}
\frac{1}{1-\frac{(1+\alpha)^2}{\alpha}s^2y} &=
\frac{\alpha}{\alpha-(1+\alpha)^2s^2y}.
\end{align*}

Solving the denominator for $\alpha$, we get two solutions:
$$\frac{1-2s^2y+\sqrt{1-4s^2y}}{2s^2y}\text{ and } \frac{1-2s^2y-\sqrt{1-4s^2y}}{2s^2y}.$$
Only the latter is a power series in $t$ with constant term $0$, which
can also be written as
$A=C(s^2y)-1$.

Thus we can apply Theorem \ref{t-lagrange1}
to get
\begin{align*}
f(y,t)=& \CT_\alpha \frac{\alpha}{\alpha-(1+\alpha)^2s^2y} (1-u^2)(1+\alpha) s^2y
\frac{1}{1+u^2(1+\alpha)^2s^2y}
 \\
=& (1-u)^2s^2y(1+A)\frac{1}{1+u^2(1+A)^2 s^2y}\frac{1}{1-2s^2y(A+1)} \\
=& (1-u^2)s^2yC(s^2y)\frac{1}{1+u^2C^2(s^2y)s^2y} \frac{1}{\sqrt{1-4s^2y}},
\end{align*}
which completes the proof.
\end{proof}

The second proof derives a different form of $F_-(y;t)$.

\begin{proof}[Second Proof of Lemma \ref{l-conj3}.]
We begin with finding the generating function of $2nb_i(n)$, which equals
 $t{\partial\over \partial t} f(y,t)$.

We claim that
\begin{align}\label{e-proof2}
\sum_{n\ge 0} {n+i\choose 2i}\frac{{4n\choose 2n}}{
{2n+2i\choose 2i}} t^{2n} =\frac{\sqrt{1+4s^2} s^{2i}}{1-4s^2},
\end{align}
where the relation between $t$ and $s$ is given in \eqref{e-of-z}.

It is easy to check that
$${n+i\choose 2i}\frac{{4n\choose 2n}}{
{2n+2i\choose 2i}} = \binom{2n-1/2}{n-i}4^{n-i}.$$

In the well-known formula
$$ \frac{C(x)^k}{\sqrt{1-4x}} =\sum_{n\ge 0} \binom{2n+k}{n} x^n,$$
by setting $x=4t^2$, and $k=2i-1/2$, we  get
$$\sum_{n\ge 0} {n+i\choose 2i}\frac{{4n\choose 2n}}{
{2n+2i\choose 2i}} t^{2n}=t^{2i}\frac{C(4t^2)^{2i-1/2}}{\sqrt{1-16t^2}}.$$
Using \eqref{e-of-z} to write the above in terms of $s$, we get
\eqref{e-proof2}.

Now we have
$$t{\partial \over \partial t} f(y;t) =
\sum_{i\ge 1}\sum_{n\ge 0} i \binom{2i}{i} {n+i\choose 2i}\frac{{4n\choose 2n}}{
{2n+2i\choose 2i}} t^{2n} y^i = \frac{2s^2y}{(1-4s^2y)^{{3/ 2}}}
\frac{\sqrt{1+4s^2}}{1-4s^2}.$$

Hence
\begin{align*}
f(y;t)=&\int \frac{2s^2y}{(1-4s^2y)^{{3/ 2}}}
\frac{\sqrt{1+4s^2}}{1-4s^2} \frac{dt}{t} \\
=& \int \frac{2s^2y}{(1-4s^2y)^{{3/ 2}}}
\frac{\sqrt{1+4s^2}}{1-4s^2} \frac{1-4s^2}{s(1+4s^2)} ds \\
=& \frac{y\sqrt{1+4s^2}}{2(1+y)\sqrt{1-4s^2y}} +\text{constant},
\end{align*}
where the constant is independent of $t$.
By setting $t=0$, and hence $s=0$, we get
$f(y;0)=\frac{y}{2(1+y)}+\text{constant}$.

Recalling equation \eqref{e-dd2}, we see that $f(y;0)=0$. Thus
the constant equals $-\frac{y}{2(1+y)}$. This gives another form of $f(y;t)$:
$$f(y;t)=
\frac{y\sqrt{1+4s^2}}{2(1+y)\sqrt{1-4s^2y}}-\frac{y}{2(1+y)}
=\frac{y(1+u^2)}{2(1+y)(1-u^2)\sqrt{1-4s^2y}}-\frac{y}{2(1+y)},$$
which is easily checked to be equal to $F_-(y;t)$ as given in \eqref{conj4}.
\end{proof}

\begin{proof}[Proof of Theorem \ref{t-conj}]
We gave a formula for the generating function
$$F(y;t)=\sum_{n\ge 0} \sum_{i\in \ZZ} a_{i,i}(2n) y^it^{2n} $$
in Lemma \ref{l-conj1}.
In Lemma \ref{l-conj2}, we showed that
$$F_-(y,t)=\sum_{n\ge 0} \sum_{i> 0} a_{-i,-i}(2n) y^it^{2n}$$
has a formula as given in \eqref{conj4}.
The proof of \eqref{conj1} is thus accomplished by Lemma \ref{l-conj3}.

For equation \eqref{conj2}, once we get
the formula \eqref{e-f+}, it is an easy exercise to show
that
$$\frac12 ((1-4s^2y)^{-1/2}-1)=\sum_{n\ge 1}\sum_{i\ge 1} 4^n\frac{i}{n}\binom{2i}{i}\binom{2n}{n-i}
y^i t^{2n}.$$

\end{proof}

{\bf Acknowledgment.}
I am very grateful to my advisor Ira Gessel, without whose help this paper
would never have been finished.
% for his helpful advices.

\end{document}